\documentclass [11pt]{article}
\usepackage{amsfonts}

\newtheorem{theorem}{Theorem}
\newtheorem{proposition}{Proposition}
\newtheorem{corollary}{Corollary}

\def\Z{{\mathbb Z}}
\def\R{{\mathbb R}}
\def\C{{\mathbb C}}

\def\Im{\mathrm{Im}\,}
\def\Re{\mathrm{Re}\,}

\begin{document}

\title{On first integrals of geodesic flows on a two-torus}
\author{I.A. Taimanov
\thanks{Sobolev Institute of Mathematics, Academician Koptyug avenue 4, 630090, Novosibirsk, Russia, and Department of Mathematics and Mechanics, Novosibirsk State University, Pirogov street 2, 630090 Novosibirsk, Russia; 
e-mail: taimanov@math.nsc.ru. \newline The work was supported by RSF (grant 14-11-00441).}}
\date{}

\maketitle

\section{Introduction}

The well-known conjecture by Kozlov (see, for instance, \cite{KD}) states that

{\sl if a geodesic flow on a two-torus admits an additional first integral polynomial in momenta, then it admits an additional first integral
 of degree $1$ or $2$ in momenta.}

Note that the geodesic flows on surfaces admitting additional first integ\-rals polynomial of
degree $k \leq 2$ in momenta were  already described in the 19th century
\cite{M,B} (see also \cite{Darboux,Birkhoff,K}).

In the recent years  many articles  appear \cite{BM11,BM15,MS,PT,S},
discussing various  approaches to this conjecture.
The conjecture is proved for the case when the conformal factor $g$
of the metric $ds^2 = g(x,y)(dx^2+dy^2)$
is a trigonometric polynomial on a torus \cite{KD}. In the general case
the conjecture stays open for every degree $k \geq 3$ of the additional integral.

We note that, by the Kozlov theorem \cite{Kozlov79}, a geodesic flow on a closed surface of genus $g>1$
does not admit an additional real analytic first integral on any positive energy level
(a generalization of this theorem for higher-dimensional manifolds was obtained in \cite{T1987,T1988}).
For a two-sphere the list of examples of integrable geodesic flow is much richer than in the torus case (see \cite{BKF}).

In the sequel, to be short we will call first integrals linear or polynomial if they
are linear or polynomial in momenta.

There are various approaches to proving this conjecture.
The equations on coefficients of the sought-for polynomials
were considered in conformal coordinates
and therewith the methods of the theory of complex-valued functions
were applied already in the 19th century.

We follow this approach in this paper which in general has a methodo\-lo\-gical character.
We will restrict a geodesic flow onto the nonzero energy level surface
for simplifying the equations.
The procedure for constructing a first integral falls into the two parts:
1) a successive formal derivation of coefficients of the polynomial;
2) a check out of closing (for this construction) reality conditions
(exactly at this step for $k=2$ the Liouville metrics are distin\-guished).

 In \S 3 it is shown that for a general metric the first part  of the
 procedure  passes through a priori only for $k \leq 7$ (Theorem 1), whereas
 for $k \geq 8$
 it is not excluded that  some conditions appear for the metrics
that allow  the successive construction of a polynomial.

The reality conditions do have different forms for the even and odd values of $k$.
In \S 4, in particular, we point out
an analogy between the reality conditions for odd $k$ and the
``dispersionless'' limits
of stationary two-dimensional soliton equations
(Remark 2 and  the equation (\ref{2hopf})).

We also consider the analogous equations for first integrals of magnetic geodesic flows.
In this case it is natural to consider  integrability only on
fixed energy levels.

Another approach to these equations, using the view point of
hydro\-dy\-na\-mic type systems, is developed in \cite{Agapov,BM12}.

 In \S 5, it is shown in particular that if the magnetic field does not
vanish or does not admit an additional linear first integral, then
 a function quadratic in momenta  can be a first integral only on one energy
level (Theorem 2).

We remark that  very few examples of integrable magnetic geodesic flows
are known. In the case when the integrals are polynomial
we have the following:

1) the flows with additional linear first integrals are completely described.
Clearly,
it is exactly the case when the conformal factor of the metric and
the magnetic field depend  on one (and the same) variable;

2) a finitely-parameterized family is known  of systems with additional quadratic first integrals \cite{DGRW};

3) it is proved that by an arbitrarily small deformation the geodesic flow
of a Liouville metric is transformed into a magnetic geodesic flow
that admits an additional quadratic integral on a given energy level;
an analytical description of these examples is unknown yet  \cite{ABM}.

We thank A.E. Mironov for helpful discussions and
the referee for helpful comments.

\section{Main notions}

Let $M$ be a Riemannian manifold homeomorphic to a two-torus.
By the uniformization theorem, it is isometric to the quotient space
$\R^2/\Gamma$, of the two-plane with a complex-valued parameter
$z = x+iy$ with respect to a lattice $\Gamma \approx \Z^2$, endowed with a conformally Euclidean metric
$$
ds^2 = g(z,\bar{z})\,(dx^2 + dy^2) = g(z,\bar{z})\,dz\,d\bar{z}.
$$

The geodesic flow on the two-torus is defined a Hamiltonian system
on the cotangent bundle $T^\ast M$ with the standard Poisson structure
in which the Poisson bracket on the coordinate functions
$x,y$, and $p_x,p_y$ (the coordinates on the fibers of the cotangent bundle, the momenta)
is as follows:
$$
\{x,p_x\} = -\{p_x,x\} = \{y,p_y\} = -\{p_y,y\} = 1
$$
and on other pairs of coordinates the Poisson bracket vanishes.
For every smooth functions $F$ and $G$, defined in a domain in $T^\ast M$, their Poisson bracket is
given by the following formula
$$
\{F,G\} = \frac{\partial F}{\partial x}\frac{\partial G}{\partial p_x} -
\frac{\partial F}{\partial p_x}\frac{\partial G}{\partial x}
+ \frac{\partial F}{\partial y}\frac{\partial G}{\partial p_y} -
\frac{\partial F}{\partial p_y}\frac{\partial G}{\partial y},
$$
which is a particular case of the general definition
$$
\{F,G\} = h^{ik} \frac{\partial F}{\partial u^i}\frac{\partial G}{\partial u^k},
$$
where $u^1,\dots,u^n$ are local coordinates on a Poisson manifold $N$ and
the skew-symmetric tensor $h^{ik}$ has the form
$$
h^{ik} = \{u^i,u^k\}, \ \ \ 1 \leq i,k \leq n.
$$

The Hamiltonian equations on a Poisson manifold $N$ are determined by a choice of a (smooth)
Hamiltonian function (the Hamiltonian)
$$
H: N \to \R
$$
and are as follows
$$
\frac{df}{dt} = \{f,H\},
$$
where the left-hand side is the evolution of an arbitrary smooth function $f$, defined in a domain in $N$,
along the trajectories of the flow. In the local coordinates
$u^1,\dots,u^n$ the equations take the form
$$
\dot{u}^i = \{u^i,H\} = h^{ik} \frac{\partial H}{\partial u^k}.
$$

For $N = T^\ast M$ the Hamiltonian
equation take the form
$$
\dot{x} = \{x,H\}= \frac{\partial H}{\partial p_x}, \ \ \
\dot{y} = \{y,H\}= \frac{\partial H}{\partial p_y},
$$
$$
\dot{p_x} = \{p_x,H\}= -\frac{\partial H}{\partial x}, \ \ \
\dot{p_y} = \{p_y,H\}= -\frac{\partial H}{\partial y}.
$$

The Hamiltonian function of the geodesic flow is equal to
$$
H = \frac{1}{2}\, h(x,y)\,(p_x^2 + p_y^2), \ \ \ h = g^{-1}.
$$
The coordinates $x$ and $y$ are not global however, since they are periodic coordinates on $M$ with respect to $\Gamma$,
the Hamiltonian equation define a dynamical system on $T^\ast M$.

Since the Poisson bracket is skew-symmetric, the Hamiltonian 
function is preserved by the flow, i.e. it is a first integral of the Hamiltonian system:
$$
\frac{dH}{dt} = \{H,H\} = 0.
$$
It determines the energy of the system.

A Hamiltonian system on a $4$-dimensional Poisson manifold $N$ is called {\it integrable}, if
there exists an additional first integral $I$:
$$
\frac{dI}{dt} = \{I,H\} = 0,
$$
which is functionally independent with $H$ almost everywhere.

More broader concept is the integrability of the fixed energy level:
a system is {\it integrable on the energy level $\{H=E\}$}, if there exists a smooth function $I$, 
defined in a neighborhood of the level surface $V = \{H=E\}$, such that
$$
\{I,H\} = 0 \ \ \ \mbox{on $V=\{H=E\}$}
$$
and $I$ and $H$ are functionally independent near $V$.

The structure of integrable hamiltonian systems (on a symplectic mani\-fold)
is described by the Liouville theorem. In particular, it states that 
regular compact submanifolds
$\{H=E, I_1=\mathrm{const},\dots,I_{k-1}=\mathrm{const}\}$ are $k$-dimen\-sio\-nal tori,
on which the Hamiltonian system is linearized (see, for instance, \cite{A,NT}).
Here the dimension on the symplectic manifold is equal to  $2k$
and for geodesic flows on a two-torus $k=2$.

These definitions are naturally generalized for all dimensions.
In this case a geodesic flow is a Hamiltonian system on the cotangent bundle
$T^\ast M$ with the Poisson bracket
\begin{equation}
\label{canon}
\{u^i,p_k\} = - \{p_k,u^i\} = \delta^i_k, \ \ \{u^i,u^k\} = \{p_i,p_k\} = 0, \ \ i,k=1,\dots,n,
\end{equation}
where $u^1,\dots,u^k$ are coordinates on the manifold, and $p_1,\dots,p_k$ are the adjoint to them coordinates in the fibers of the cotangent bundle.
We recall that coordinates in a domain of a symplectic manifold are call canonical if
for them the Poisson bracket takes the form (\ref{canon}).

The Hamiltonian function of the geodesic flow (of a Riemannian metric) is defined as
$$
H(u,p) = \frac{1}{2}g^{ij}(u)p_i p_j,
$$
where $ds^2 = g_{ij}(u)du^idu^j$ is the Riemannian metric and $g^{ij}g_{jm} = \delta^i_m$ for $i,m=1,\dots,k$.

For the geodesic flow on any Riemannian manifold $M$ the integrability on a nonzero energy level, in fact, implies 
the integ\-ra\-bi\-lity, because the restrictions of the flow onto different nonzero energy levels are orbitally isomorphic.
Therewith for a study of topological obstructions to the integrability it is more convenient to work with the energy level 
$V$ which is compact if the configuration space $M$ is compact \cite{T1987}.

Everywhere in the sequel we will assume that the energy level $H=E$ under consideration corresponds to a positive value of energy:
$$
E > 0,
$$
i.e. it is not empty as for $E<0$ and does not consist of all stationary points of the flow as for $E=0$.

\section{Equations on first integrals of geodesic flows}

The Hamiltonian $H(u,p)$ of a geodesic flow is homogeneous in momenta
in canonical coordinates on the cotangent bundle.
Therefore we may apply to the geodesic flow the following proposition which is checked
by straightforward computations.

\begin{proposition}
\label{pro1}
If, in canonical coordinates, $F(u,p)$ is a polynomial of $q$-th degree in momenta $p_1,\dots,p_n$ and 
the Hamiltonian function $H(u,p)$ is a polynomial of $p$-th degree in momenta, then $\{F,H\}$ is a polynomial of $(p+q-1)$-th degree in momenta.
\end{proposition}

\begin{corollary}
\label{cor1}
If $F(u,p)$ is a real analytic first integral of the geodesic flow and
 $F(u,p) = \sum_{k=0}^\infty F_k(u,p)$ is its decomposition into a series,
where
 $F_k(u,p)$ is a polynomial of $k$-th degree in momenta, $k=0,1,\dots$, then
 all functions $F_k(u,p)$ are first integrals of the geodesic flow.
\end{corollary}

Let us return to geodesic flows on the two-torus $M = T^2$
with the conformal parameter
$z=x+iy$.

We  extend the space of smooth functions on $T^\ast M$  to the space of
smooth complex-valued functions on the same phase space.
For a pair of functions $F$ and $G$ their Poisson bracket is defined
as in the real-valued case by the formula
$$
\{F,G\} = \frac{\partial F}{\partial x}\frac{\partial G}{\partial p_x} -
\frac{\partial F}{\partial p_x}\frac{\partial G}{\partial x}
+ \frac{\partial F}{\partial y}\frac{\partial G}{\partial p_y} -
\frac{\partial F}{\partial p_y}\frac{\partial G}{\partial y}.
$$
Since the fibers of the cotangent bundle $T^\ast M$ are diffeomorphic to the complex line $\C$, 
the complex-valued functions on the fibers it is naturally to consider as
functions of complex-valued parameters
$$
p_z = \frac{1}{2}(p_x - ip_y), \ \ \ p_{\bar{z}} = \frac{1}{2}(p_x + ip_y).
$$
The Poisson brackets of complex-valued parameters are equal to
$$
\{z,p_z\} = \{\bar{z},p_{\bar{z}}\} = 1, \ \ \ \{z,p_{\bar{z}}\} = \{\bar{z},p_z\}=0.
$$
The Hamiltonian $H(u,p)$ of the geodesic flow of a metric
$ds^2 = g(x,y)\,(dx^2+dy^2) =  g(z,\bar{z})\, dz\,d\bar{z}$
takes the form
$$
H(u,p) = \frac{1}{2}h(x,y)(p_x^2+p_y^2) = 2 h(z,\bar{z})p_z p _{\bar{z}}, \ \ h=g^{-1}.
$$

By Corollary \ref{cor1}, the existence of real analytic first integrals of geodesic flows is equivalent
to the existence of the first integrals that are homogeneous polynomials
in momenta. Hence when we are looking for first integrals
we restrict ourselves by the first integrals of such a form.

A real-valued function $f$ on $T^\ast M$, which is a homogeneous polynomial
of $k$-th degree in momenta, has the form
\begin{equation}
\label{pol}
f = a_{k,0} p_z^k + a_{k-1,1}p_z^{k-1}p_{\bar{z}} + \dots + a_{1,k-1} p_z p_{\bar{z}}^{k-1} +
a_{0,k}p_{\bar{z}}^k,
\end{equation}
where
$$
a_{k-l,l} = a_{k-l,l}(z,\bar{z}), \ \ a_{k-l,l} = \bar{a}_{l,k-l}, \ \ l=0,\dots,k.
$$

The following proposition  is checked by straightforward computations:

\begin{proposition}
A function $f$ of the form (\ref{pol}) is a first integral
of the geodesic flow on the energy level
$$
E = 2 h p_z p_{\bar{z}} = \mathrm{const},
$$
if and only if the function $\{f,H\}$, which is  polynomial in $p_z$
and $p_{\bar{z}}$, vanishes after the substitution
\begin{equation}
\label{subst}
p_z p_{\bar{z}} = \frac{E}{2h}.
\end{equation}
\end{proposition}

We remark that after the substitution (\ref{subst})  $f$ takes
the polynomial form
\begin{equation}
\label{pol2}
F = a_k p_z^k + \dots + a_1 p_z + a_0 + \bar{a}_1 p_{\bar{z}} + \dots + \bar{a}_k p_{\bar{z}}^k, \ \ a_0=\bar{a}_0,
\end{equation}
where $a_i=a_i(z,\bar{z}), i=0,\dots,k$.
The following proposition is evident.

\begin{proposition}
A function $f$ of the form (\ref{pol})
is the first integral of the geodesic flow on the energy level $H=E$ if and only if the function
$F$ obtained from $f$ by applying the substitution (\ref{subst}) is the first integral of the flow
on the same energy level.
\end{proposition}

Let us compute the results of the successive application to polynomials of the form
(\ref{pol2}) of the Poisson bracket with $H$ and the substitution (\ref{subst}). We denote this operation by
$$
F \to \{F,H\}_E
$$
and note that it sends polynomials of the form (\ref{pol2}) into
polynomials of the same form:

1) $F = a(z,\bar{z}) p_z^m, m>0$:
$$
\{a p_z^m, 2h p_z p_{\bar{z}}\}_E = a\{p_z^m, 2hp_z p_{\bar{z}}\}_E + p_z^m \{a, 2hp_z p_{\bar{z}}\}_E =
$$
$$
= am p_z^{m-1}\{p_z, 2h p_zp_{\bar{z}}\}_E + p_z^m \frac{\partial a}{\partial z}\{z, 2h p_z p_{\bar{z}}\}_E +
p_z^m \frac{\partial a}{\partial \bar{z}}\{\bar{z}, 2h p_z p_{\bar{z}}\}_E =
$$
$$
-ma\frac{E}{h}\frac{\partial h}{\partial z} p_z^{m-1} +E \frac{\partial a}{\partial z}p_z^{m-1} + 2h \frac{\partial a}{\partial \bar{z}} p_z^{m+1};
$$

2) $F = \bar{a}(z,\bar{z}) p_{\bar{z}}^m, m>0$:
$$
\{ \bar{a} p_{\bar{z}}^m, 2h p_z p_{\bar{z}}\} = \overline{\{a p_z^m, 2h p_z p_{\bar{z}}\}}_E =
$$
$$
= -m\bar{a}\frac{E}{h}\frac{\partial h}{\partial \bar{z}} p_{\bar{z}}^{m-1} +
E \frac{\partial \bar{a}}{\partial \bar{z}}p_{\bar{z}}^{m-1} + 2h \frac{\partial\bar{a}}{\partial z} p_{\bar{z}}^{m+1};
$$

3) $\bar{F} = a(z,\bar{z}), a =\bar{a}$:
$$
\{a, 2hp_z p_{\bar{z}}\}_E = 2h p_z \{a,p_{\bar{z}}\} + 2h p_{\bar{z}} \{a,p_z\} =
$$
$$
=
2h \frac{\partial a}{\partial \bar{z}} p_z + 2h \frac{\partial a}{\partial z}p_{\bar{z}}.
$$

These formulas imply the following proposition.

\begin{proposition}
\label{proint}
1) A function $F_k$ of the form (\ref{pol2}) is a first integral of the geodesic flow on 
the energy level $H=E$ if and only if
$\{F_k,H\}_E=0$, where
$$
\{F_k,H\}_E =
\sum_{m=1}^{k+1} p_z^m \left( 2h \frac{\partial a_{m-1}}{\partial \bar{z}} - 
(m+1) a_{m+1} \frac{E}{h} \frac{\partial h}{\partial z} + E \frac{\partial a_{m+1}}{\partial z}\right) +
$$
$$
\left(E \frac{\partial a_1}{\partial z} - a_1\frac{E}{h}\frac{\partial h}{\partial z} +
E\frac{\partial \bar{a}_1}{\partial \bar{z}} - \bar{a}_1 \frac{E}{h}\frac{\partial h}{\partial \bar{z}}\right) +
$$
$$
+
\sum_{m=1}^{k+1} p_{\bar{z}}^m \left( 2h \frac{\partial \bar{a}_{m-1}}{\partial z} - 
(m+1) \bar{a}_{m+1} \frac{E}{h} \frac{\partial h}{\partial \bar{z}} + E \frac{\partial \bar{a}_{m+1}}{\partial \bar{z}}\right).
$$

2) If a function $F_k$ of the form (\ref{pol2}) is a first integral of the geodesic flow on the energy level
$H=E$, then both components $F_{\mathrm{even}}$ and $F_{\mathrm{odd}}$ of its decomposition
$$
F_k = F_{\mathrm{even}} + F_{\mathrm{odd}}
$$
 into the sum of functions, of the same form, which are invariant and anti-invariant with respect to the transformation
$(p_z,p_{\bar{z}}) \to (-p_z,-p_{\bar{z}})$ are first integrals of the geodesic flow on this energy level.
The functions $F_{\mathrm{even}}$ are $F_{\mathrm{odd}}$ are polynomials of even and odd degrees in momenta, respectively.
\end{proposition}

\begin{corollary}
The coefficient $a_k$ of the decomposition of the first integral $F_k$,
at the highest degree of $p_z$ is holomorphic in $z$:
\begin{equation}
\label{kolokol}
\frac{\partial a_k}{\partial \bar{z}}=0.
\end{equation}
\end{corollary}

{\sc Remark 1.}
For $k=2$ the equation (\ref{kolokol}) was pointed out already in the 19th century: see, for instance,
\cite[p. 42]{Darboux}, where the local problem was considered and,
by using (\ref{kolokol}),
a quadratic first integral was reduced to the form $a_2=1$ by
a local change of coordinates. It was substantially used by Birkhoff  \cite{Birkhoff} and in \cite{K},
where, in particular, the nonexistence of additional first integrals,
real analytic in momenta on closed surfaces of genus $g>1$ was derived from
the nonexistence of not-everywhere-vanishing holomorphic $k$-differentials $a(z)(dz)^k$.

Application of Proposition \ref{proint} to a search for metrics
on tori with additional first integrals polynomial in momenta
meets the two difficulties:

a) the formal solvability of the equations on $a_i, 0 \leq i \leq k$,

b) the fulfillment of the ``closing'' reality conditions for derived solutions.

The appearing closing conditions are different in form for even and odd values of
$k$. Let us elaborate that.

\subsection{$k=2q+1$}

By part 2 of Proposition \ref{proint}, it is enough to consider the case when
a polynomial $F_k$ has nonzero coefficients only at the odd degrees of $p_z$
and $p_{\bar{z}}$.

The system of equations on the coefficients of the polynomial
$$
F_k = a_{2q+1} p_z^{2q+1} + a_{2q-1} p_z^{2q-1} + \dots + a_1 p_z + \bar{a}_1 p_{\bar{z}} + \dots +
\bar{a}_{2q-1} p_{\bar{z}}^{2q-1} + \bar{a}_{2q+1} p_{\bar{z}}^{2q+1}
$$
is written as a condition of vanishing of coefficients at $p_z^{2m}, 0 \leq m \leq q+1$,
of the decomposition of $\{F_k, H\}_E$. In this case we divide every coefficient by $h$
and, after replacing $h$ by $g=h^{-1}$ everywhere, we obtain
\begin{equation}
\label{intodd0}
\frac{\partial a_{2q+1}}{\partial \bar{z}} = 0, \ \ \ m=q+1.
\end{equation}
\begin{equation}
\label{intodd}
\frac{\partial a_{2m-1}}{\partial \bar{z}} + \frac{E}{2}\left((2m+1) a_{2m+1}\frac{\partial g}{\partial z} +
g \frac{\partial a_{2m+1}}{\partial z}\right) = 0, \ \ \ 1 \leq m \leq q;
\end{equation}
\begin{equation}
\label{intodd2}
\Re \frac{\partial (g\,a_1)}{\partial z} = 0, \ \ \ \ m=0.
\end{equation}

Since, by (\ref{kolokol}), $a_k = \mathrm{const}$ as a holomorphic double-periodic function, the value of $a_{2q+1}$
can be mapped to any nonzero constant by
a linear change of a variable $z \to z/c$.

Then we haves to successively solve the equations for $a_{2m-1}$, decreasing the value of $m$ from $q$ to $1$.
In general the question of the successive solvability of this system is open.

The necessary and sufficient condition for the solvability of the equation on $a_{l}$, when the values of $a_{l+2},\dots,a_{k}$ are given,
is as follows.
Let us write down the equation (\ref{intodd}), which corresponds to this case, in the general form:
\begin{equation}
\label{solv}
\frac{\partial a_n}{\partial \bar{z}}  = - \frac{E}{2}\left((n+2) a_{n+2}\frac{\partial g}{\partial z} + g \frac{\partial a_{n+2}}{\partial z}\right).
\end{equation}
We denote by $C^\infty(M)$ the space of smooth functions on the torus $M=\R^2/\Gamma$ and by
$C_0^\infty(M)$ its subspace which consists of all functions $f$ such that
$$
\int_M f \, dx\,dy = 0.
$$
This condition is equivalent to the vanishing of the free term in the Fourier decomposition of $f$.
On $C_0^\infty(M)$ the operators $\partial = \frac{\partial}{\partial z}$ and
$\bar{\partial} = \frac{\partial}{\partial \bar{z}}$ are invertible up to constants.
For simplicity, we demonstrate that for the case of the square lattice $\Gamma = 2\pi \Z^2$.
The Fourier decomposition of  $f$ is as follows:
$$
f(x,y) = f_0 + \sum_{-\infty \leq k,l \leq \infty} f_{k,l} e^{i(kx+ly)},
$$
and exactly when $f_0 = 0$ the functions
$$
\partial^{-1}f(x,y) = \sum_{-\infty \leq k,l \leq \infty}  \frac{2 f_{k,l}}{ik+l}  e^{i(kx+ly)}
$$
and
$$
\bar{\partial}^{-1}f(x,y) = \sum_{-\infty \leq k,l \leq \infty}  \frac{2 f_{k,l}}{ik-l}  e^{i(kx+ly)}
$$
satisfy the equations
\begin{equation}
\label{solv2}
\partial u = f, \ \ \ \bar{\partial} v = f,
\end{equation}
respectively. All other solutions of the equations from $C_0^\infty(M)$ differ from these by constants.
It is clear that for
$f_0 \neq 0$ the equations are unsolvable in $C^\infty(M)$.

Let us return to the equations (\ref{solv}). We already mentioned that $a_k=\mathrm{const}$.
The equation (\ref{solv2}) is solvable as an equation on $a_n$ with an assumption
that  $a_{n+2},\dots,a_k$
are known if and only if its right-hand side
$$
(n+2) a_{n+2}\frac{\partial g}{\partial z} + g \frac{\partial a_{n+2}}{\partial z} =
(n+1) a_{n+2}\frac{\partial g}{\partial z} + \frac{\partial (g\,a_{n+2})}{\partial z}
$$
(for brevity, we take it up to a multiple constant)
lies in $C_0^\infty(M)$. This is equivalent to the equality
\begin{equation}
\label{condition1}
a_{n+2}\frac{\partial g}{\partial z} \in C_0^\infty(M).
\end{equation}
The last condition is satisfied for $n=k-2, k-4, k-6$. Indeed:

1) for $a_k = \mathrm{const}$ and $n=k-2$ the conditions reduce to
$\frac{\partial g}{\partial z} \in C_0^\infty(M)$, which is evident, and
the equation on $a_{k-2}$
takes the simple form
$$
\frac{\partial a_{k-2}}{\partial \bar{z}} = -\frac{\lambda}{2} \frac{\partial g}{\partial z},
$$
where, for brevity, we assume
$$
\lambda = k\,a_k\,E;
$$

2) for $n=k-4$ we have
$$
a_{k-2}\, \frac{\partial g}{\partial z} = - \frac{2}{\lambda}\, a_{k-2}\, \frac{\partial a_{k-2}}{\partial \bar{z}} =
-\frac{1}{\lambda}  \frac{\partial a_{k-2}^2}{\partial \bar{z}}
$$
and the right-hand side lies in $C_0^\infty(M)$ as a derivative of a smooth function on the torus;

3) for $n=k-6$
$$
a_{k-4}\, \frac{\partial g}{\partial z} = - \frac{2}{\lambda}\, a_{k-4}\, \frac{\partial a_{k-2}}{\partial \bar{z}} =
-\frac{2}{\lambda} \frac{\partial (a_{k-4}a_{k-2})}{\partial \bar{z}} + \frac{2}{\lambda}a_{k-2}
\frac{\partial a_{k-4}}{\partial \bar{z}}.
$$
From this computation it is easy to notice that the condition (\ref{condition1}) is satisfied if and only if
\begin{equation}
\label{condition2}
a_{k-2} \frac{\partial a_{n+2}}{\partial \bar{z}} \in C_0^\infty(M)
\end{equation}
and for checking the last condition we may replace $\frac{\partial a_{n+2}}{\partial \bar{z}}$
by the right-hand side of the equation (\ref{solv}) on $a_{n+2}$. For $n=k-6$ we get
$$
-\frac{2}{E}\,a_{k-2} \frac{\partial a_{k-4}}{\partial \bar{z}} =
(k-2)a_{k-2}\frac{\partial g}{\partial z} + g \frac{\partial a_{k-2}}{\partial z}
= (k-2)a_{k-2}^2 \frac{\partial g}{\partial z} + \frac{1}{2}\frac{\partial a_{k-2}^2}{\partial z} =
$$
$$
\left(k-\frac{5}{2}\right) a_{k-2}^2\,\frac{\partial g}{\partial z} + \frac{1}{2} \frac{\partial (g\,a_{k-2}^2)}{\partial z} =
- \frac{2k-5}{\lambda}a_{k-2}^2\frac{\partial a_{k-2}}{\partial \bar{z}} + \frac{1}{2} \frac{\partial (g\,a_{k-2}^2)}{\partial z} =
$$
$$=
- \frac{2k-5}{3\lambda} \frac{\partial a_{k-2}^3}{\partial \bar{z}} + \frac{1}{2} \frac{\partial (g\,a_{k-2}^2)}{\partial z} \, \in \, C_0^\infty{M}.
$$

For $n=k-8$ the condition (\ref{condition2}) takes the form
$$
a_{k-2} \frac{\partial a_{k-6}}{\partial \bar{z}} = -\frac{E}{2} 
\left( (k-4)a_{k-2}\,a_{k-4}\,\frac{\partial g}{\partial z} + a_{k-2} g \frac{\partial a_{k-4}}{\partial z}\right)
$$
and we do not know how to show that the right-hand side of it lies in $C_0^\infty(M)$.
It is possible that this and the subsequent equations on $a_n$ are solvable
not for all metrics $g$. We proved the following

\begin{theorem}
\label{thsolv}
For every metric $ds^2 = g(x,y)(dx^2+dy^2)$ on a two-torus the system of equations
(\ref{intodd0}) and (\ref{intodd}) is successively solvable for
$a_k=\mathrm{const}, a_{k-2}$, $a_{k-4}$, and $a_{k-6}$.

On every step a solution $a_n$ is obtained up to a constant.
\end{theorem}

\begin{corollary}
The equations $\{F_k,H\}_E=0$ on $F_k$ are solvable for all metrics for $k\leq 7$.
\end{corollary}

These theorem and corollary are valid for all (even and odd) values of $k$, because the equations have the same form. From
(\ref{condition1}) it follows that the solvability of the subsequent equations on $a_n$ does not depend on the choices of $a_m, m > n$, up to constants.

A polynomial $F_k$, which satisfies (\ref{intodd0}) and (\ref{intodd}), gives an additional first integral if the reality condition
(\ref{intodd2}) is satisfied. It would not be superfluous to write down this condition one more time:
$$
\Re \frac{\partial (g\,a_1)}{\partial z} = 0.
$$

Up to recently all known cases of additional polynomial integrals of odd degrees
are reduced to the case, then
$g(x,y) = g(y), \ \ \ F_k = \mathrm{const}\cdot p_x^k$.
We consider it in detail in \S \ref{subsec4.1}.
For $k=1$ the reality condition (\ref{intodd2}) is evidently satisfied, because $a_1= E \,g(y)$.

\subsection{$k=2q$}

It follows from Proposition \ref{proint} that the equations on $F_k$ are obtained as the vanishing conditions for
the coefficients at  $p_z^{2m+1}, 0 \leq p$, in the decomposition of $\{F_k,H\}_E$:
$$
\frac{\partial a_{2p}}{\partial \bar{z}} = 0, \ \ \ \ m=p;
$$
$$
\frac{\partial a_{2m}}{\partial \bar{z}} + \frac{E}{2}\left((2m+2) a_{2m+2}\frac{\partial g}{\partial z} +
g \frac{\partial a_{2m+2}}{\partial z}\right) = 0, \ \ \ 0 \leq m \leq p-1.
$$
Solutions of this system are derived successively starting with $a_k = \mathrm{const}$.
Theorem \ref{thsolv} gives a sufficient condition ($k \leq 6$) for the solvability of the system for all metrics $g$.

In difference with the case of odd $k$ the reality condition appears not as the vanishing of one of coefficients of the decomposition of
$\{F,H\}_E$, but as the reality condition for the polynomial $F_k$:
\begin{equation}
\label{evenreal}
a_0 = \bar{a}_0.
\end{equation}

It is nontrivial already for $k=2$ and in this case all solutions of the system are given by the Liouville metrics (see \S \ref{subsec4.2}).
For arbitrary even $k=2q$ all known cases of additional integrals of even degrees relate to the Liouville metrics
and $F_{2q} = \mathrm{const}\cdot F_2^q$.

When we solve the equations on $F_k$ we obtain in general case a polynomial
$F_k$, whose coefficients $a_m$ are polynomials in $E$. Therewith if in the final formula we substitute
$E$ by $E = H = 2hp_zp_{\bar{z}}$, we obtain a homogeneous polynomial of degree
$k$ in $p_z$ $p_{\bar{z}}$, which is a first integral for all values of energy.
We demonstrate this evident fact below for the case of the Liouville metrics
(\S \ref{subsec4.2}).

\section{First integrals of low degrees}

First, we expose the well-known results on geodesic flows with additional
linear and quadratic integrals.

\subsection{$k=1$}
\label{subsec4.1}

By Proposition \ref{proint}, a first integral takes the form
$$
F_1 = ap_z + \bar{a}p_{\bar{z}}
$$
and therewith
$$
\frac{\partial a}{\partial {\bar{z}}} = 0,\ \ \
E \frac{\partial a}{\partial z} - a\frac{E}{h}\frac{\partial h}{\partial z} +
E\frac{\partial \bar{a}}{\partial \bar{z}} - \bar{a} \frac{E}{h}\frac{\partial h}{\partial \bar{z}} = 0.
$$
The first equation implies that
$a = \mathrm{const}$
and by change of a variable $z \to z/a$ we may reduce the polynomial to the form $a=1$.

For $a=1$ the second equation is written as
$$
\frac{\partial h}{\partial x} = 0, \ \ \ g=h^{-1},
$$
i.e. the metric takes the form
$$
ds^2  = g(y)(dx^2 + dy^2),
$$
and the first integral is a canonical momentum corresponding to the coordinate $x$:
$F_1 = p_x$.

\subsection{$k=2$}
\label{subsec4.2}

It is enough to consider the case
$$
F_2 = a_2 p_z^2 + a_0 + \bar{a}_2 p_{\bar{z}}^2, \ \ a_0=\bar{a}_0,
$$
because, if $F_{\mathrm{odd}}\neq 0$, then there exists a linear first integral and this case was already considered.
Proposition \ref{proint} implies that
$$
\frac{\partial a_2}{\partial {\bar{z}}} = 0,
$$
$$
2h \frac{\partial a_0}{\partial \bar{z}} - 2a_2 \frac{E}{h}\frac{\partial h}{\partial z} +
E\frac{\partial a_2}{\partial z}=0.
$$
Again we conclude that $a_2=\mathrm{const}$ and by a change of a variable
$z \to z/\sqrt{a_2}$ we reduce the polynomial to the case $a_2=1$.
Then the second equation reads
$$
\frac{\partial a_0}{\partial \bar{z}} - \frac{E}{h^2}\frac{\partial h}{\partial z} = 0,
$$
which is rewritten in terms of the metric $g=h^{-1}$ as follows
\begin{equation}
\label{e2}
\frac{\partial a_0}{\partial \bar{z}} + E \frac{\partial g}{\partial z} = 0.
\end{equation}
From the last equation the function $a_0$ is obtained by the inversion of the operator
$\bar{\partial}=\frac{\partial}{\partial \bar{z}}$ uniquely up to constant. However an additional reality condition as to be satisfied:
$$
a_0 = \bar{a}_0.
$$
This condition distinguishes the class of the metrics admitting an additional quadratic integral, so-called
{\it Liouville metrics}. To find it we differentiate the left-hand side of (\ref{e2}) in $z$
and obtain
$$
\frac{\partial^2 a_0}{\partial z \, \partial \bar{z}} + E \frac{\partial^2 g}{\partial z^2} = 0.
$$
A real-valued function $a_0$ meets this equation if and only if
$$
\Im \frac{\partial^2 g}{\partial z^2} = \frac{1}{2} \frac{\partial^2 g}{\partial x \, \partial y} = 0
$$
(here we keep in mind that $g$ is also real-valued).
By the D'Alembert formula, $g(x,y)$ has the form
$$
g(x,y) = v(x) + w(y),
$$
$a_0$ is found (up to a summand $\mathrm{const}\cdot E$)
in the form
$$
a_0 = -E \left(\frac{\partial^2}{\partial z \, \partial \bar{z}}\right)^{-1} \frac{\partial^2 g}{\partial z^2} =
-E \Delta^{-1} (v^{\prime\prime}(x) - w^{\prime\prime}(y)) = -E(v(x) - w(y)),
$$
and we obtain a first integral
$$
F_2 = p_z^2 + p_{\bar{z}}^2 - E (v(x) - w(y)).
$$
We substitute into this formula
$$
E = \frac{1}{2g} (p_x^2+p_y^2) = \frac{p_x^2+p_y^2}{2(v+w)}
$$
and derive that
$$
F = \frac{1}{2}(p_x^2-p_y^2)  - \frac{(p_x^2+p_y^2)(v-w)}{2(v+w)} =
\frac{p_x^2\, w(y) - p_y^2 \, v(x)}{v(x)+w(y)}.
$$
This is an additional quadratic first integral of the geodesic flow of the Liouville metric
$$
ds^2 = (v(x)+w(y))(dx^2 + dy^2).
$$
We note that although we  looked for a first integral on the nonzero energy level $H=E$,
the energy enters linearly in the final formulas and, by sub\-sti\-tu\-ting its value by $H$, we obtain a first integral for all energy levels.
The non-uniqueness in he choice of $a_0$ consists in adding an additional term $\mathrm{const}\cdot E$, which finally results in adding to $F$
a term of the form  $\mathrm{const}\cdot H$.

In the particular case $v(x)=0$ we obtain a metric whose geodesic flow admits a linear first integral $F_1=p_x$ and $F_2 = F_1^2$.

We notice that the period lattice $\Gamma$ is not necessarily rectangular.
It is only necessary that the projection $\R^2 \to \R^2/\Gamma$
send the lines $x = \mathrm{const}$ and $y=\mathrm{const}$ on the plane $\R^2$
to closed curves \cite{BN}.

 \subsection{$k=3$}

The equations for coefficients reduce to the equation on $a_1$ and therewith we can, by using, if necessary, a linear change of $z$, put $a_3$ 
to be equal to any nonzero constant.
It is convenient to put $a_3=\frac{1}{3}$ and we obtain the system of two equations:
$$
\frac{\partial a_1}{\partial \bar{z}} + \frac{E}{2} \frac{\partial g}{\partial z}=0,
$$
$$
\Re \frac{\partial (g\,a_1)}{\partial z} = 0.
$$

{\sc Remark 2.}
Let us consider this system of equation from another point of view.
To make it compatible with notations usual for soliton theory we put
$$
u = \frac{E\,g}{2}, \ \ \ v = -a_1,
$$
and rewrite it as
\begin{equation}
\label{2hopf}
\frac{\partial (uv)}{\partial z} + \frac{\partial (u\bar{v})}{\partial \bar{z}} = 0, \ \ \ \
\frac{\partial v}{\partial \bar{z}} =
\frac{\partial u}{\partial z}.
\end{equation}
In soliton theory it is known the Novikov--Veselov equation (NV) \cite{NV}
$$
\frac{\partial u}{\partial t} = \frac{\partial^3 u}{\partial z^3} + \frac{\partial^3 u}{\partial \bar{z}^3} +
\frac{\partial (uv)}{\partial z} + \frac{\partial (u\bar{v})}{\partial \bar{z}},
$$
where $v$ is determined by the equation
$$
\frac{\partial v}{\partial \bar{z}} =
\frac{\partial u}{\partial z},
$$
which is uniquely solvable in the class of fast decaying functions on
$\R^2$ and in the class on functions with vanishing mean value on a two-torus.
The Novikov--Veselov equation has two natural one-dimensional reductions:

\begin{enumerate}
\item
$u=u(x), v = u$: in this case it reduces to the Korteweg--de Vries equation (KdV)
$$
\frac{\partial u}{\partial t} = \frac{1}{4}\frac{\partial^3 u}{\partial x^3} + 2 u \frac{\partial u}{\partial x}.
$$
Therefore the NV equation is a two-dimensional generalization, of the KdV equation, which differs from the Kadomtsev--Petviashvili equation;

\item
$u=u(y), v=-u$: in this case the right-hand side of the NV equation is trivial
and any function $u(y)$, which depends only on $y$, is a stationary solution of the NV equation
(this is valid also for all equation from the NV hierarchy which under the reduction of the first type
becomes the KdV hierarchy.
\end{enumerate}

\noindent
In the ``dispersionless'' limit the KdV equation
$$
t \to \varepsilon^{-1} t, \ \ x \to \varepsilon^{-1} x, \ \ \varepsilon \to 0
$$
the KdV equation becomes the Hopf equation
$$
\frac{\partial u}{\partial t} = 2 u \frac{\partial u}{\partial x},
$$
which, as it is known, has no non-stationary smooth solutions different from constants: $u = \mathrm{const}$.
The system  (\ref{2hopf}) may be considered as the equation for stationary solutions of 
the two-dimensional generalization of the Hopf equation, i.e. of the ``dispersionless'' limit of the NV equation.
Solutions which depends only on $y$ exist and in terms of integrable geodesic flows they correspond to metrics of the form
$g(y)(dx^2+dy^2)$. Therewith the first integral of third degree is proportional to $F_1^3 = p_x^3$.
It could be possible that only constant solutions of this equation are smooth.
This analogy is an argument for the following conjecture

{\sl the equation (\ref{2hopf}) has no smooth nontrivial (depending substantially on $x$) solutions.}

\noindent
Similar analogies with the Hopf equation appear in the cases of other odd degrees $k >3$.

{\sc Remark 3.} In  \cite{BM11,BM15}, by using semigeodesic coordinates on a torus, it is shown that the condition for 
the existence of an additional polynomial first integral, i.e. the equations on coefficients of a polynomial in momenta, 
are written as systems of hydrodynamic type which also appear in soliton theory.

\subsection{$k=4$}

Put, for simplicity,  $a_4=\frac{1}{4}$. The final system takes the form
$$
\frac{\partial a_2}{\partial \bar{z}} + \frac{E}{2}\frac{\partial g}{\partial z} = 0,
$$
$$
\frac{\partial a_0}{\partial \bar{z}} + \frac{E}{2}\left(2a_2 \frac{\partial g}{\partial z} +
g \frac{\partial a_2}{\partial z}\right) = 0,
$$
with $a_0=\bar{a}_0$.

In difference with the case of odd degrees $k$ we can not draw analogies between these equations and soliton equations and
introduce arguments pros and cons the existence of solutions which differ from ones given by the Liouville metrics.

\subsection{Some remarks}

The integrability of the geodesic flows of Liouville metrics was established in \cite{L}.
The complete classification of (two-dimensional) metrics admitting locally additional linear or quadratic first intergal was
obtained by Massieu \cite{M}. By (\ref{kolokol}), locally the coefficient at the highest degree can be transformed to constant
by a change of variables and therefore the equations for $F_k$ have the same form as in Proposition \ref{proint},
and for low degrees $k=1,2$ the Massieu classification reduces to metrics depending on one variable and
to the Liouville metrics. For $k\geq 3$ such a classification is not obtained. For $k \leq 2$
the classification of such geodesic flows on a two-torus was apparently first time derived in  \cite{K}.

In \cite{L} Liouville considered the geodesic flow of a metric
$g\,dz\,d\bar{z}$ on the energy level
$$
H = 2g^{-1}p_z p_{\bar{z}} = E
$$
as a motion of a particle in the Euclidean metric and the potential field $U(x,y) = E\,g(x,y)$ on the zero energy level:
$$
H -E = 0 \longrightarrow (H-E)g = \frac{1}{2}p_zp_{\bar{z}} - Eg = 0.
$$
Clearly this is one of the forms of the Maupertuis principle.

For the systems describing a motion in a potential field $U(x,y)$ on a two-torus with the Euclidean metric
$ds^2= dx^2 + dy^2$ there is a conjecture which states that if there is an additional polynomial first integral 
(for all energy levels), then there exists such on integral of degree $k \leq 2$. For $k\leq 2$ such systems are 
completely described: for $k=1$ the potential $U$ depends on one variable: $U=U(x)$,
and for $k=2$ it has the ``Liouville''form: $U(x,y) = V(x)+W(y)$. Up to recently this conjecture is proved for low degrees:
$k \leq 5$ \cite{Byaly,DK,DKT,Mironov10}.

\section{Magnetic geodesic flows}

The inclusion of a magnetic filed into a system consist in
the deformation of the Poisson bracket.  In the canonical coordinates $x,y,p_x$, and $p_y$
the inclusion of a magnetic field results in the replacement of the relation $\{p_x,p_y\} = 0$ by
$$
\{p_x,p_y\} = B
$$
where the Poisson brackets for other pairs of coordinates are preserved \cite{Novikov1982}.
Here $B dx \wedge dy$ is a  $2$форма on a surface $M$.
In the multidimensional case the Poisson brackets of pairs of momenta
take the form
$$
\{p_i,p_k\} = B_{ik},
$$
where $B = \sum_{i<k} B_{ik} du^i \wedge du^k$ is a closed $2$-form on a manifold $M$ (in the two-dimensional case
the form $B\,dx\wedge dy$ is always closed). The form $B$ describes the magnetic field.

For complex-valued coordinates on a surface the Poisson bracket takes the form
$$
\{z,p_z\}=\{\bar{z},p_{\bar{z}}\} = 1, \ \ \ \{z,\bar{z}\}=0, \ \ \ \{p_z,p_{\bar{z}}\} = iB.
$$
The Hamiltonian of the magnetic geodesic flow with respect to this Poisson bracket is the same as for the geodesic flow:
$$
H = 2 h p_z p_{\bar{z}}, \ \ h^{-1}=g, \ \ ds^2 = g(z,\bar{z})\,dz\,d\bar{z}.
$$

A magnetic field may be exact or non-exact: in the first case there exists a $1$-form $A=A_i \,du^i$, on the manifold,
such that $B=dA$, in the second case there exists a smooth mapping $\varphi: N \to M$ 
of a two-dimensional oriented closed manifold $N$ into $M$ such that $\int_N \varphi^\ast (B) \neq 0$. If $M$ is a closed oriented surface, then
$B$ is non-exact if and only if $\int_M B \neq 0$.

Since for magnetic geodesic flows the Poisson brackets of momenta are nontrivial,
for such systems the statements of Proposition
\ref{pro1}, Corollary \ref{cor1}, and part 2 of Proposition \ref{proint} do not hold.
The simplest example of violations of the analogs of these statements is given by the flat metric
$ds^2 = dx^2+dy^2$ and the constant magnetic field $B\, dx\wedge dy, B=\mathrm{const}\neq 0$ on the two-torus
$M = \R^2/2\pi \Z^2$. In this case an additional first integral is equal to
$$
F = \cos \left(\frac{p_x}{B}-y\right),
$$
it is real analytic and its components homogeneous in momenta are not first integrals of the flow.
The trajectories of this flow on every nonzero energy level are circles of the same radius
and $F$ is the value of a smooth function (on the two-torus) at the center of the circle.

The analogous example for close hyperbolic surfaces is more complicated:
if $ds^2$ is a metric of constant curvature $K=-1$, $d\mu$ is the area form of this metric, and
$B dx \wedge dy = Kd\mu$, then for $E>\frac{1}{2}$ the magnetic geodesic flow is chaotic, and for
$E<\frac{1}{2}$ it is integrable, its trajectories on the universal covering are the geodesic circles of constant radii
and for the first integral we may take the value of a fixed smooth function (on the compact surface) at the center of the circle
\cite{T2004}.

Let us consider the problem of finding magnetic geodesic flows
which on nonzero energy levels $H=E\neq 0$ have additional first integrals of the form
\begin{equation}
\label{mint}
F_k = a_k p_z^k + \dots + a_1 p_z + a_0 + \bar{a}_1 p_{\bar{z}} + \dots + \bar{a}_k p_{\bar{z}}^k.
\end{equation}
The decomposition $\{F_k,H\}_E$ into the powers of momenta is as follows
\footnote{As it was pointed out to us by the referee, the formula (\ref{expansion}) and Proposition \ref{exactness}
which follows from it were obtained in  \cite{Bolotin}.}
\begin{equation}
\label{expansion}
2h \frac{\partial a_k}{\partial \bar{z}} p_z^{k+1} + 2h \left(\frac{\partial a_{k-1}}{\partial \bar{z}} +
ik\,a_k\,B\right)p_z^k + \dots.
\end{equation}

In the sequel, we will consider the case when $M = \R^2/\Gamma$ is a two-torus with
a complex-valued parameter $z$.
As in the case of geodesic flows we conclude it follows from
$\{F_k,H\}_E=0$ that
$$
a_k  = \mathrm{const},
$$
and the vanishing of the coefficient at $p_z^k$ implies the equality
$$
\frac{i}{k\,a_k} \frac{\partial a_{k-1}}{\partial \bar{z}} = B, \ \ B=\bar{B},
$$
from which we derive that
$$
\frac{\partial \alpha}{\partial x} = \frac{\partial \beta}{\partial y}, \ \ \ 
-\frac{1}{2\,k\,a_k}\left(\frac{\partial \alpha}{\partial y}+\frac{\partial \beta}{\partial x}\right) = B,
$$
where $a_{k-1}=\alpha+i\beta, \alpha=\bar{\alpha}, \beta=\bar{\beta}$.
This immediately implies

\begin{proposition}
\label{exactness}
If a magnetic geodesic flow on a two-torus on some nonzero energy level
$\{H=E\neq 0\}$ admits a polynomial first integral $F_k$ (of the form (\ref{mint})), then the magnetic field is exact:
$$
B \, dx\wedge dy = d\left(\frac{1}{k\,a_k}(\alpha\,dx - \beta\,dy)\right).
$$
\end{proposition}

We have also the following elementary statement which is implied by
the equation $\{F_1,H\}_E=0$.

\begin{proposition}
Let a magnetic geodesic flow on a two-torus admit on some nonzero energy level
$\{H=E \neq 0\}$ a linear first integral $F_1$. Then in convenient conformal coordinates we have
$$
F_1 = p_z + a_0 + p_{\bar{z}},  \ \ \ ds^2 = g(z,\bar{z})\,dz\,d\bar{z}, \ \
\frac{\partial a_0}{\partial x} = \frac{\partial g}{\partial x}=0, \ \ B = - \frac{1}{2}\frac{\partial a_0}{\partial  y}.
$$
Moreover $F_1$ is a first integral on all energy levels.
\end{proposition}

Specific examples of magnetic geodesic flows with linear (in momenta) first integrals were almost not studied
however the possibility of their explicit description allows in some cases to find interesting dynamical properties
 \cite{T15}.

For $k=2$ the situation is more interesting.

\begin{proposition}
Let a magnetic geodesic flow on a two-torus admits on some nonzero energy level
$\{H=E=\neq 0\}$ an additional quadratic first integral $F_2$.
Then in convenient conformal coordinates
$$
F_2 = a_2 p_z^2 + a_1 p_z + a_0 + \bar{a}_1 p_{\bar{z}} + \bar{a}_2 p_{\bar{z}}^2,  \ \ \ ds^2 = g(z,\bar{z})\,dz\,d\bar{z},
$$
where $a_0, a_1, B$, and $g$ satisfy the equations
$$
\frac{\partial a_1}{\partial \bar{z}} + 2i\,a_2 B = 0,
$$
\begin{equation}
\label{mint2}
\frac{\partial a_0}{\partial \bar{z}} + E\,a_2\, \frac{\partial g}{\partial z} + 2ia_1 B = 0,
\end{equation}
$$
\frac{\partial (g\, a_1)}{\partial z} + \frac{\partial (g\,\bar{a}_1)}{\partial \bar{z}} = 0,
$$
where $a_0 = \bar{a}_0, B = \bar{B}, g=\bar{g}$, and $a_2=\mathrm{const}$.
\end{proposition}

The system (\ref{mint2}) immediately follows from Proposition \ref{proint} and the addi\-tio\-nal nontrivial contribution
of the Poisson bracket on momenta: $\{p_z,p_{\bar{z}}\}=iB$. Therewith this additional contribution of the magnetic field
does not affect the highest term in $p_z$ in the decomposition of $\{F_k,H\}_E$ and the equation (\ref{kolokol})
is satisfied.

We apply this Proposition for deriving some facts.

\begin{theorem}
\label{th2}
For a magnetic geodesic flow with a real analytic magnetic field on a two-torus a quadratic in momenta
 function of the form
$$
F_2 =  a_2 p_z^2 + a_1 p_z + a_{01} + a_{02}H + \bar{a}_1 p_{\bar{z}} + \bar{a}_2 p_{\bar{z}}^2,
$$
can be a first integral of the flow on two different nonzero energy levels if the magnetic field vanishes everywhere
of the magnetic geodesic flow admits a linear first integral.
\end{theorem}

{\sc Proof.}
Since the integrability on one energy level implies  $a_2=\mathrm{const}$, by a change of variable we reduce
$F$ to the case $a_2=1$.
Let us assume that the flow is integrable on the different energy levels $E_1$ and $E_2$.
By (\ref{mint2}), we have
$$
\frac{\partial (a_{01} + a_{02}E)}{\partial \bar{z}} + E_n \, \frac{\partial g}{\partial z} + 2ia_1 B = 0, \ \ \ n=1,2,
$$
and, since $E_1 \neq E_2$, this implies two equalities:
$$
\frac{\partial a_{01}}{\partial \bar{z}} + \frac{\partial g}{\partial z} = 0, \ \ \
\frac{\partial a_{02}}{\partial \bar{z}} + 2ia_1 B = 0.
$$
However, by (\ref{mint2}),
\begin{equation}
\label{magnet}
B = \bar{B} = \frac{i}{2}\frac{\partial a_1}{\partial \bar{z}} = -\frac{1}{4}\left(\frac{\partial \alpha}{\partial y} + \frac{\partial \beta}{\partial x}\right)
+\frac{i}{4}\left(\frac{\partial \alpha}{\partial x} - \frac{\partial \beta}{\partial y}\right),
\end{equation}
where $\alpha = \Re a_1, \beta = \Im a_1$,
and we conclude that
$$
\frac{\partial}{\partial \bar{z}} \left(a_{02} - \frac{1}{2}a_1^2\right) = 0,
$$
which implies
$$
a_{02} - \frac{1}{2}a_1^2 = \mathrm{const}.
$$
Since $a_{02} = \bar{a}_{02}$, we have $\Im a_1^2 = \alpha \beta = \mathrm{const}$.
Let us consider two cases:

1) $\alpha \beta = C = \mathrm{const} \neq 0$. Then $\beta = C/\alpha$, by (\ref{magnet}),
$$
\frac{\partial \alpha}{\partial x} = \frac{\partial \beta}{\partial y} = - \frac{C}{\alpha^2}\frac{\partial \alpha}{\partial y}
$$
from which we conclude that
$$
\frac{\partial \alpha}{\partial y} = -\frac{\alpha^2}{C}\frac{\partial \alpha}{\partial x}.
$$
This equation as the Hopf equation does not admit nonconstant smooth solutions.
Therefore $\alpha$ and $\beta$ are both constant and $B=0$.

2) $\alpha\beta=0$. In this case either $\alpha=0$, either $\beta=0$  everywhere
(exactly where we need the real analyticity of $B$ and, therefore, of $a_1$).
For $\alpha = 0$, by (\ref{magnet}), we have
$$
\frac{\partial \beta}{\partial y} = 0, \ \ \ B=B(x), \ \ \ a_1 = a_1(x).
$$
For $\beta=0$ analogously we derive that $B=B(y)$ and $g=g(y)$.
In both cases the magnetic geodesic flow admits a linear first integral.
Theorem is proved.

The proof of the following statement is completely analogous.

\begin{theorem}
If a Liouville metric $g = v(x) + w(y)$ depends substantially on both variables, then for every nontrivial (i.e. not vanishing everywhere) magnetic field $B$ on every nonzero energy level the magnetic geodesic field does not admit
a first integral of the form
$$
F_2 = p_z^2 + a_1 p_z + a_0 + \bar{a}_1 p_{\bar{z}} + p_{\bar{z}}^2.
$$
\end{theorem}

{\sc Proof.} Let us assume that for some nontrivial energy level such a first integral exists.
By (\ref{mint2}), we have
$$
\frac{\partial a_0}{\partial \bar{z}} + E\, \frac{\partial g}{\partial z} + 2ia_1 B = 0.
$$
By differentiating the left-hand side in $z$ and substituting $B= \frac{i}{2}\frac{\partial a_1}{\partial \bar{z}}$,
we obtain
$$
\frac{\partial^2}{\partial z\,\partial \bar{z}} \left(a_0 - \frac{a_1^2}{2}\right) = - E \frac{\partial^2 g}{\partial z^2}.
$$
Since $g$ is a Liouville metric, the right-hand side of the last equation
is real-valued as well as  $a_0$. Therefore,
$\Im a_1^2 = \mathrm{const}$ and, by applying the same reasonings as in the proof of Theorem \ref{th2},
we arrive at contradiction. Theorem is proved.

Until recently neither local, no global (on closed surfaces) classification of magnetic geodesic flows
admitting on a fixed energy level an additional quadratic first integral are not obtained.

It is known only one explicit example which is derived, by using the Maupertuis principle,
from an integrable case of a motion of a particle in magnetic and potential fields \cite{DGRW},
which we expose below for the completeness of an exposition.

Recently in \cite{ABM} it was shown that every geodesic flow of a Liouville metric by an arbitrarily small deformation is transformed into a magnetic geodesic flow which admits a fixed energy level an additional quadratic first integral.

{\sc Example} (\cite{DGRW}). We expose it in the spirit of the equations (\ref{mint2}). Put $a_2=1$.
Let the magnetic field have the ``Liouville'' form:
$$
B = B_1(x)+B_2(y).
$$
since for the existence of an additional first integral it is necessary that the magnetic field is exact,
i.e. its mean value over the torus vanishes:
$\int B\,dx\wedge dy=0$, the functions $B_1$ are $B_2$ are the derivatives of periodic functions.
In this case, without loss of generality, we may assume that
$$
B = v^{\prime\prime}(x) + w^{\prime\prime}(y).
$$
From (\ref{mint2}) the vector potential of the magnetic field is easily derived in the form
$$
a_1 = -4(w^\prime(y) + i v^\prime(x))
$$
and the condition $\Re \frac{\partial (g\,a_1)}{\partial z}=0$ reduces to the equation
$$
w^\prime \frac{\partial g}{\partial x} + v^\prime \frac{\partial g}{\partial y}
 =0,
$$
whose solutions are as follows
$$
g(x,y) = f(v(x)-w(y)),
$$
where $f(u)$ is an arbitrary differentiable function of one variable.
It is left to find $a_0$ from the equation
\begin{equation}
\label{eq}
\frac{\partial a_0}{\partial \bar{z}} + E \frac{\partial g}{\partial z} - \frac{1}{2}\frac{\partial a_1^2}{\partial \bar{z}}= 0
\end{equation}
and therewith the reality condition $a_0=\bar{a}_0$ has to be satisfied. Up to a constant $a_0$ is uniquely
determined by the last equation as well as by the equation which is obtained by differentiating
both sides of (\ref{eq}) in $z$:
$$
\frac{\partial^2 a_0}{\partial z\,\partial\bar{z}} + E \frac{\partial^2 g}{\partial z^2} - \frac{1}{2}\frac{\partial^2 a_1^2}{\partial z\, \partial\bar{z}} = 0.
$$
The last equation admits a real-valued solution $a_0$ if and only if
$$
\Im \left(E \frac{\partial^2 g}{\partial z^2} - \frac{1}{2}\frac{\partial^2 a_1^2}{\partial z\, \partial\bar{z}}\right)=0.
$$
Let us substitute into this condition the expressions for $a_1$ and $g$ and obtain
$$
Ef^{\prime\prime}v^\prime w^\prime - 8\Delta (v^\prime w^\prime) =
Ef^{\prime\prime}v^\prime w^\prime - 8(v^{\prime\prime\prime}w^\prime+ v^\prime w^{\prime\prime\prime}) = 0.
$$
For
$$
f(u) = \alpha_3 u^3 + \alpha_2 u^2 + \alpha_1u+\alpha_0, \ \ \ \alpha_0,\dots,\alpha_3 \in \R,
$$
the last equality is written as
\begin{equation}
\label{dgrw1}
E(6\alpha_3 (v-w)+2\alpha_2)v^\prime w^\prime - 8 (v^{\prime\prime\prime}w^\prime+ v^\prime w^{\prime\prime\prime}) = 0.
\end{equation}
If
$$
6E \alpha_3 vv^\prime + 2\gamma v^\prime  - 8 v^{\prime\prime\prime}= 0,
$$
$$
-6E\alpha_3 ww^\prime + 2(\alpha_2-\gamma) w^\prime - 8 w^{\prime\prime\prime} = 0,
$$
then (\ref{dgrw1}) holds. Note that the elliptic function $\varphi$, defined by the formula
$$
\varphi^{\prime\,2} = c_3 \varphi^3 + c_2\varphi^2 + c_1 \varphi + c_0,
$$
where $c_0,\dots,c_3$ are constants,
satisfies the equation
$$
\varphi^{\prime\prime\prime} = \frac{3}{2}c_3 \varphi \varphi^\prime + c_2 \varphi^\prime.
$$
We obtain the final answer in the form
$$
g(x,y) = \alpha_3(v(x)-w(y))^3 + \alpha_2 (v(x)-w(y))^2 + \alpha_1 (v(x)-w(y)) +\alpha_0,
$$
$$
v^{\prime\, 2} = \frac{\alpha_3E}{2} v^3 + \frac{\gamma}{4} v^2 + \gamma_1 v + \gamma_0,
$$
$$
w^{\prime\,2} = - \frac{\alpha_3E}{2}w^3 + \left(\frac{\alpha_2-\gamma}{4}\right)w^2 + \delta_1 w + \delta_0,
$$
where $\alpha_0,\dots,\alpha_3,\gamma,\gamma_0,\gamma_1,\delta_0,\delta_1$ are real-valued constants.
For brevity, we skip the expression for $a_0$.

{\sc Remark 4.} The general integrability problem of magnetic geodesic flows
has been studied much  less than its analogue for geodesic flows.
Mainly physically interesting low-dimensional examples were considered.
For exam\-ple, the paper \cite{DGRW} concerns natural mechanical systems on the Euclidean plane with
magnetic and potential fields and the integrable case, which is found in it,
is naturally equivalent to a magnetic geodesic flow on a two-torus (see above).
In \cite{Bolotin}, in particular, there was established the non-integrability on fixed energy levels of natural
mechanical systems on two-manifolds in the presence of gyroscopic and potential fields with sufficiently many
singularities of the Newton type.
More general configuration spaces recently began to be considered:
in \cite{Efimov1,Efimov2} the complete integrability of magnetic geodesic flows on compact 
simply-connected homogeneous symplectic manifolds with magnetic fields given by the symplectic forms 
was established, and additional algebraic reasonings together with the deformation of the momentum map 
introduced in \cite{Efimov1,Efimov2} allow to substantially expand the class of integrable examples \cite{BJ0,BJ}.

\end{document}